\documentclass[12pt]{amsart}
\usepackage{amsfonts}
\usepackage{amssymb}

\advance\textwidth2.3cm
\advance\oddsidemargin-1.7cm
\advance\evensidemargin-1.7cm
\advance\textheight 2.6cm
\advance\topmargin-1.5cm

\newskip\smallh

\def\aside#1\par{\par
\setbox0\hbox{#1}
\smallh0.6\hsize
\relax
\ifnum\wd0>\smallh\else\smallh\wd0\fi
\rightline{
\vrule$\,$\vrule\  \vbox{\rightskip0cm plus 1cm
 \hsize\smallh\parindent 0cm #1}\ 
\vrule$\,$\vrule}\par}

  \renewcommand{\c}{\mathfrak c}
  
  \renewcommand{\r}{\mathfrak r}
  \renewcommand{\d}{\mathfrak d}
 
\newcommand{\s}{\mathfrak s}

\renewcommand{\b}{\mathfrak b}
\newcommand{\e}{\mathfrak e}

\newcommand{\measure}{\mathbb L}       
\newcommand{\category}{\mathbb K}

\newcommand{\F}{{\mathcal F}}

\newcommand{\R}{\mathfrak R}
  \renewcommand{\S}{\mathfrak S}

\newcommand{\cov}{\operatorname{cov}}

\newcommand{\on}{\mathord{\restriction}}
\newcommand{\dom}{\operatorname{dom}}

\newcommand{\aut}{\operatorname{Aut}}

\newcommand{\sym}{\operatorname{Sym}}
\newcommand{\supp}{\operatorname{supp}}
\newcommand{\fin}{\text{fin}}

\newcommand{\forces}{\Vdash}

\def\pre#1^#2{{}^{#2}{#1}}
\def\oo{[\omega]^\omega}

\def\ccc{\pre 2^\omega }  

\newtheorem{thm}{Theorem}[section]
\newtheorem{theorem}[thm]{Theorem}
\newtheorem{Theorem}[thm]{Theorem}
\newtheorem{lemma}[thm]{Lemma}
\newtheorem{Lemma}[thm]{Lemma}

\theoremstyle{definition}
\newtheorem{definition}[thm]{Definition}
\newtheorem{Definition}[thm]{Definition}

\newtheorem{Fact}[thm]{Fact}

\newtheorem{corollary}[thm]{Corollary}

\newtheorem{Claim}[thm]{Claim}
\newtheorem{Conclusion}[thm]{Conclusion}
\newtheorem{Problem}[thm]{Problem}
\newtheorem{Remark}[thm]{Remark}

\def\itm#1 {\item[(#1)]}

\newcommand{\Cal}[1]{{\mathcal #1}}
\newcommand{\om}{\omega}
\newcommand{\su}{\subseteq}

\long\def\ignore#1\endignore{}

\title{Rules and Reals}
\author{Martin Goldstern}
\thanks{The first author is supported by the Austrian Science
  Foundation (FWF)}
\author{Menachem Kojman}
\thanks{The second author
  was partially supported by an NSF grant no. 9622579}
\date{June 1997}

\begin{document}
\begin{abstract}

A  ``$k$-rule" is a sequence $\vec A=((A_n,B_n):n<\om)$ of pairwise disjoint
sets $B_n$, each of cardinality $\le k$ and subsets $A_n\su B_n$. A
subset $X\su \om$ (a ``real'') follows a rule $\vec A$ if for infinitely
many $n\in \om$, $X\cap B_n=A_n$.

There are obvious cardinal invariants resulting from this definition:
the least number of reals needed to follow  all $k$-rules, $\s_k$,  and the
least number
of $k$-rules without a real following  all of them, $\r_k$.

Call $\vec A$ a {\em bounded} rule if $\vec A$ is a $k$-rule for some
$k$. Let $\r_\infty$ be the least cardinality of a set of bounded rules
with no  real  following all rules in the set.

We prove the following:
$\r_\infty\ge\max(\cov(\category),\cov(\measure))$ and $\r=\r_1\ge \r_2=\r_k$
for all $k\ge 2$. However, in the Laver model, $\r_2<\b=\r_1$. 

An application of $\r_\infty$ is in Section 3: we show that below $\r_\infty$
one can find proper extensions of dense independent families which preserve a
pre-assigned group of automorphisms. The original motivation for discovering
rules was an attempt to construct a maximal homogeneous family
over $\om$. The consistency of such a family is still open.

\end{abstract}

\maketitle

%
%
\section*{Introduction}

In the present paper we present new cardinal invariants which resulted
from investigations of homogeneous families. These numbers have
intrinsic interest (in fact we regard it as surprising that those
numbers have not been discovered earlier).

In Section 1 we discuss cardinal invariants related to
``$k$-rules.''
A  {\em $k$-rule}
 is a sequence $\vec A=((A_n,B_n):n<\om)$ of pairwise disjoint
sets $B_n$, each of cardinality $\le k$, and subsets $A_n\su B_n$. A
subset $X\su \om$ (a ``real'') {\em follows}  a rule $\vec A$ if for infinitely
many $n\in \om$, $X\cap B_n=A_n$.

A rule $\vec A$ is {\em bounded} if it is a $k$-rule
for some $k\in \om$.

The obvious cardinal invariants related to rules are the following: the
least number of reals needed to follow all $k$-rules, $\s_k$, and the
least number of $k$-rules with no real following all of them,
$\r_k$. Let $\r_\infty$ be the least number of bounded rules with no real
following all of them. 

We compare the $\r_k$s and $\r_\infty$ among themselves and to well known
cardinal invariants: covering of category, covering of Lebesgue measure,
$\r$, $\b$, $\d$ 
and the evasion numbers $\e_k$ which were studied by Blass and Brendle. We
prove: 
\begin{itemize}
\itm a
 $\max(cov(\category), cov(\measure)) \le \r_\infty$;
\itm b
$\r=\r_1\ge
\r_2=\r_k$ for all $k\ge 2$;
\itm c
  $\s_2\le\e_2$;
\itm d
 $\r_\infty\le
\min(\r_2,\d)$. 
\end{itemize}

In Section 2 we prove the consistency of $\r_2<\b$.

In Section 3 we show that below $\r_\infty$ one can properly extend an
independent family of subsets of $\omega$ preserving a prescribed group
of automorphisms. This is the relevance of $\r_\infty$ to the behavior of
homogeneous families under inclusion, which was the original motivation
for the discovery of rules.
 

%
%
%
%
\section{Rules}

\begin{definition}
\begin{itemize}
\itm 1 A {\em rule} is a sequence $\vec A = 
        ( A_n, B_n: n \in \omega)$, where the sets
        $B_n$ are disjoint and finite, and for all $n$, $A_n \subseteq
        B_n\subseteq\om$. 
\itm 2
We say that $X\in \oo $ {\em follows} the rule $\vec A$ if there are
        infinitely many $n$ with $X \cap B_n = A_n$; otherwise $X$ is
        said to {\em avoid} $\vec A$. 
\itm 3
For $k \in \omega$ we say that $\vec A$ is a $k$-rule if all sets
        $B_n$ have size $\le k$.   We say that $\vec A$ is a
        {\em bounded} rule if $\vec A$ is a $k$-rule for some
        $k$.    
\item 4
More generally, for any function
        $f:\omega \to \omega$ we say that $\vec A$ is an $f$-rule if
        for all $n$, $|B_n| \le f(n)$. 
        We say that $f$ is a ``slow'' function if 
        $$ \sum_{n=0}^\infty 2^{-f(n)} = \infty,$$
and we say that $\vec A$ is a slow rule if it is an $f$-rule for some 
slow $f$.
\end{itemize}
\end{definition}

\begin{definition}
\begin{enumerate}
\itm 1 For $k \in \omega $ let $\r_k:= 
\min \{ |\R|: $ there is no $X$ which follows all $k$-rules from $\R$ $\}$.
        (Similarly $\r_f$, when $f:\omega \to \omega$.) 
\itm 2 Dually, let $\s_k:= 
\min \{ |\S|: $ every $k$-rule is followed by some $X\in \S$  $\}$.
\itm 3 We let $\r_\infty = \min
\{ |\R|: $ there is no $X$ which follows all bounded rules from
$\R$\nobreak$\}$.
\end{enumerate}
\end{definition}

We remark that $2^k$ trivially bounds the least cardinality of
a set of $k$-rules with the property that every real follows some rule
in the set.

Recall that the ``splitting'' number $\s$
and the ``reaping'' number $\r$ are defined as follows:

\begin{definition}
If $s, X \in \oo $, then we say that $s$ ``splits'' $X$ if
$s$ divides $X$ into two infinite parts, i.e., $s\cap X$ and
$(\omega-s)\cap X$ are both infinite. 

\begin{itemize}
\itm 1 $\s:= \min \{|\S|: \S \subseteq \oo $,
         every $X\in \oo $ is split by some
        $s\in \S$ $\}$
\itm 2 $\r:= \min \{|\R|$: $\R \subseteq \oo $, there is no
        $X\in \oo $ which 
        splits all $r\in \R$ $\}$
\end{itemize}
\end{definition}

\begin{Fact}
$\cdots \le \r_3 \le \r_2 \le \r_1 = \r$, 
        and $\s \le \s_2 \le \s_3 \le \cdots$.   However, 
$\s_1 = 2$, witnessed by $\S =\{\emptyset, \omega \}$. 
\end{Fact}

\begin{Theorem}
\begin{enumerate}
\itm  a Let $( N,{\in})$ 
be a model of ZFC* (a large enough fragment of ZFC). 
If a real $X$ follows all rules from $N$, then $X$ is Cohen over
        $N$. (Conversely, a Cohen real over $N$ follows all rules from $N$.)
\itm b 
        If $X$ is random over $N$, then $X$ follows all slow
        rules from $N$ (so in particular, all  bounded rules). 

\itm c $\max(cov(\category), cov(\measure)) \le \r_\infty$.
($\cov(\category)$ is the smallest number of first category
sets needed to cover the real line.  $\cov(\measure )$ is
defined similarly using measure zero sets.)
\end{enumerate}
\end{Theorem}

\begin{proof}

(a): Assume that $X \subseteq \omega$ follows all rules from
$N$.  We claim that $\chi_X$, the characteristic function of
$X$, is a Cohen real over $N$, that is, the set $\{ \chi_X\on
n: n \in \omega\}$ is generic for the forcing notion
$\pre2^{<\omega}$.

To verify this claim, consider any 
nowhere dense tree
 $T \subseteq \pre 2^{<\omega}$ in $N$. 
We have to check that $\chi_X$ is not a branch of $T$. 

Using the fact that $T$ is nowhere dense (and $T$ is in $N$)
we can by induction (in $N$!)\  
find sequences $( n_i: i < \omega)$ and $( \eta_i : i <
\omega)$ such that for all $i< \omega$ we have: 
\begin{enumerate}
\item $n_i < n_{i+1} $, $\eta_i \in \pre 2^{[n_i, n_{i+1})} $
\item For all $\nu \in \pre 2^{n_i}$, $\nu \cup \eta_i \notin T$. 
\end{enumerate}

Now let $B_i:= [n_i, n_{i+1})$, $A_i =
\{ k: \eta_i(k)=1 \}$. Our
assumption tells us 
 that $X$ follows the rule $( A_i, B_i: i \in
\omega)$.  So for some $i$ we have $X \cap B_i = A_i$, and hence
$\chi_X \supseteq \eta_i$. Hence $\chi_X$ is not a branch of
$T$.

This concludes the proof of (a).

The converse to (a)  is obvious.

(b) is also easy: Let $X_n:= \{ X: X \cap B_n \not= A_n \}$.
For $n\not=m$, the sets $X_n$ and $X_m$ are  independent (in
the probabilistic sense), and $\mu(X_n) = 1-2^{-f(n)}$, where
$\mu$ is the Lebesgue measure on ${\Cal P} (\omega) \simeq \pre
2^\omega $.  Hence $\mu(\bigcap_{n>m} X_n) = \prod_{n>m}
(1-2^{-f(n)}) = 0$. 

(c) follows from (a) and (b). 

\end{proof}

\begin{Theorem}[Shelah]\label{rkr2}
For $k\ge 2$, $\r_k = \r_2$ (and similarly, $\s_2 = \s_k$). 
\end{Theorem}

\begin{proof}
We will show that $\r_k=\r_{k+1}$: Let $N_0$ be sufficiently closed 
(say, a model of ZFC*, but closed under some recursive functions is
        sufficient)  of 
size $<\r_{k}$;  we have to show that there is a real that follows
all $k+1$-rules from $N_0$. 

We define a sequence $(N_i, C_i: i \le k)$ such that 
$N_i \cup \{ C_i \} \subseteq N_{i+1}$, each $N_i$ is sufficiently
        closed and of the same cardinality as $N_0$, and $C_i$ follows
        all $k$-rules from $N_i$.  

Now let $C$ be the ``average'' of the $C_i$: 
        $m \in C$ iff $m$ is in ``most'' of the $C_i$'s, or formally:
$$ C:= \{ m \in \omega: | \{i\le k: m \in C_i \} | > (k+1)/2 \}$$

Now we check that $C$ indeed follows all $k+1$-rules from $N_0$. 

Let $(A_n, B_n: n \in \omega)$ be a $k+1$-rule in $N_0$. 
For $0\le i \le k$ we let $(A^i_n, B^i_n: n \in \omega)$ be the 
$k$-rule obtained by removing the each $i$th element of $B_n$.   
That is, letting $\{b_n^0, \ldots, b_n^k\}$ be the increasing enumeration
of $B_n$ we let $B_n^i:= B_n \setminus \{b^i_n\}$, 
$A_n^i:= A_n \cap B_n^i$.

Let $E_0:= \omega$. 
For $0\le i \le  k$ let 
$$E_{i+1}:= \{ n \in E_i:  
                B^{i}_n \cap C_i = A^{i}_n\}, $$
i.e., $E_{i+1}$ is the set of indices on which $C_i$ follows                 
 the rule $(A_n^i, B_n^i: n \in E_i)$.    Note that $E_i \in N_i$ 
and $C_{i} \in N_{i+1}$.    By the choice of $C_i$ we know that 
each $E_{i+1}$ is infinite. 

We conclude the proof by showing that for $n \in E_{k+1}$ we
have $A_n = B_n \cap C$.    Let $n \in E_{k+1}$ 
(so also $n \in E_i$ for all $i \le k$), and $m \in B_n$.   
Say $m = b_n^j$.   Then for $i \not= j$ we have  $m \in B_n^i$, 
so $m \in A_n^i  \Leftrightarrow m \in C_i$. 

Hence the cardinality of the set  
$ \{ i \le k: m \in C_i \}$ is either in 
$\{0,1\} $ (iff $m\notin A_n$) or in $\{ k, k+1\}$.  
In any case we 
 get $m \in C $ iff $m \in A_n$.  
 So $A_n = B_n \cap C$. 

\end{proof}

\begin{theorem}
$\r_\infty \ge \min(\r_2, \d)$.  In particular, if $\r_2 \le
\d$ then $\r_\infty = \r_2$.
\end{theorem}
\begin{proof}
Let $N$ be a model of ZFC* of cardinality $< \min(\r_2, \d)$.   We
will find a real $X$ following all bounded rules from $N$.

Define sequences $(N_i:i<\omega)$, $(X_i:i<\omega)$ 
satisfying the following conditions:
\begin{enumerate}
\item $N_0 = N$.                 
\item  $N_i$ is a model of ZFC*, $N_{i-1} \cup\{X_i\}
                \subseteq N_i $.
\item $|N_i| = |N_0|$.                 
\item $X_i$ follows all $i$-rules (and hence also all
$j$-rules for $j\le i$) from $N_{i-1}$. 
\end{enumerate}
Let $N_\omega$ be a model of size $|N_0|$ containing
$(N_i:i<\omega)$ and $(X_i:i<\omega)$.  Since $|N_\omega | <
\d$ we can find a strictly
increasing  function $f$ that is not dominated by any
function from $N_ \omega $.

Define $X \subseteq  \omega $ by requiring $X \cap  (   
f(i-1), f(i)]  = X_i \cap  (f(i-1), f(i)]$.
We claim that $X$ follows all
bounded rules from $N$.

To complete the proof, consider an arbitrary $k$-rule
$(A_n,B_n:n \in  \omega )$  from
$N$.   We may assume $\min \bigcup_n B_n > f(k)$.
We define sequences $(E_i:k \le i < \omega ) $
satisfying the following conditions
for all $i \ge k$. 
\begin{enumerate}
\item $\forall n \in E_i$ $B_n \cap X_i = A_n$.
\item $E_i \in N_i$.
\item $E_{i+1} \subseteq E_i$. 
\end{enumerate}
We can carry out this construction, because
$(A_n, B_n: n \in E_i)$ is a rule in  $N_i$, so we just choose
$E_{i+1}$ to witness that $X_{i+1}$ follows this rule. 

Now let $n_i:= \min E_i$.  Clearly the function $i \mapsto n_i$
is in $N_\omega $.   So we can find infinitely many $j$ such that $f(j) >
\max B_{n_j} $.

We claim that for each such $j$,
 $X \cap B_{n_j} = A_{n_j}$. For all  $i \in [k,j]$
 we have $n_j \in E_i$, so $X_i \cap B_j  = A_j$.  Note that
$B_j  \subseteq [f(k), f(j))$, so we also have $X \cap
B_j=A_j$. 

\end{proof}

\begin{Problem}  Is $\r_\infty < \r_2$ consistent?
\end{Problem}

We remark that in the random real model we have
$\r_2=\cov(\measure) = \c = \r_\infty$,
$\d=\aleph_1$. So one cannot hope to prove $\r_2\le \d$.

\bigskip

We now consider the invariant that is dual to $\r_k$, and we
compare it with the well-known ``evasion'' number.

\begin{definition}
$(\pi,D)$ is a $k$-predictor, if $D$ is an infinite
subset of $ \omega $, $\pi = (\pi_n:n \in D)$, 
$\pi_n $ a function from $\pre k^n$ to $k$. 

We say that $f\in \pre k^\omega$ evades $(\pi, D)$ if there are
infinitely many $\ell\in D$ such that $f(\ell) \not= \pi_\ell(f\on \ell)$. 

$$ \e_k := \min \{ |N|: \forall \pi\,  \exists f\in N: 
        \hbox {$f$ evades $\pi$} \}$$
\end{definition}

Brendle in \cite{evpr} investigated these and other cardinal
invariants and showed  that all $\e_k$ are equal to each
other. 

The following construction connects rules with predictors. 

\begin{definition}
Let $R = (A_n, B_n: n \in \omega)$ be a 2-rule.  Define a
$2$-predictor $(\pi_R, D_R)$ as follows:
\begin{enumerate}
\item $D_R = \{\max B_n: n \in \omega \}$
\item If $\ell=\max B_n$, and $|A_n|=1$, then $\pi_\ell(f)=f(\min
B_n)$ for all $f \in \pre 2^\ell$.  Otherwise, $\pi_\ell(f)=1-f(\min
B_n)$. 
\end{enumerate}
\end{definition}

\begin{lemma}
Let $X \subseteq \omega$.   If $\chi_X$ 
 evades $\pi_R$, then either $X$ or $\omega \setminus X$
follows $R$. 
\end{lemma}
\begin{proof}
Let $\ell_n:= \max B_n$, $i_n = \min(B_n)$ for all $n$.

$X$ evades $\pi_R$, so there are infinitely many $n$ such that 
$X(\ell_n) \not= \pi_{\ell_n}(X\on \ell_n)$.

\subsubsection*{Case 1:} There are
infinitely 
many such $n$ where in addition $|A_n| = 1$. 

So for each such $n$, $X(\ell_n) \not= \pi_{\ell_n}(X\on \ell_n) = 
X(i_n)$.  So $X(\ell_n) \not= X(i_n)$, so $X \cap B_n$ must be either 
$A_n$ or $B_n\setminus A_n$.     One of the two alternatives holds
infinitely often.  Hence, either there are infinitely many $n$
such that $X\cap B_n = A_n$, or there are infinitely many $n$ such
        that $(\omega \setminus X) \cap B_n = A_n$. 

\subsubsection*{Case 2:} There are
infinitely many such $n$ with 
$X(\ell_n) \not= \pi_{\ell_n}(X\on \ell_n)$,
where in addition $|A_n| = 2$,  i.e., $A_n = B_n$. 
So for each such $n$, $X(\ell_n) \not= \pi_{\ell_n}(X\on \ell_n) = 
1-X(i_n)$.  So $X(\ell_n) = X(i_n)$, so $X \cap B_n$ must be either 
$B_n$ or $\emptyset$.     One of the two alternatives holds
infinitely often.  So again we either get infinitely many $n$
such that $X\cap B_n = A_n$, or  infinitely many $n$ such
        that $(\omega \setminus X) \cap B_n = A_n$. 

\subsubsection*{Case 3:}  For infinitely many $n$ as above
we have $A_n=\emptyset$.  Similar to the above.

\end{proof}

\begin{corollary}
$\s_2 \le \e_2$
\end{corollary}

\begin{proof}
Let $N$ be a model (of set theory) 
witnessing $\e_2$, i.e., for every 2-predictor 
$\pi$ there is a function $f\in N$ evading $\pi$.    

Let $R$ be any 2-rule.   There is $X\in N$ evading $\pi_R$, so either
$X$ or $\omega \setminus X$  (both in $N$) follows $R$. 
\end{proof}

\begin{Remark} $\s \le \e_2$ is known. 
Brendle showed that $\s<\s_2$ is consistent (unpublished).
\end{Remark}

%
%
%
%
\section{Consistency of $\r_2<\r$}

We show here in contrast to theorem \ref{rkr2}
that $\r$ is not provably equal to $\r_2$.  Moreover, whereas
$\b\le \r$ is provable in ZFC (see \cite{Vaughan} for a
collection of results
on cardinal invariants), we show that $\r_2 < \b$ is
consistent with ZFC.

The following definition is standard: 

\begin{Definition}
\begin{enumerate}
\item  $S$ is a {\em slalom} iff $\dom(S)= \omega$ and for all
        $n\in \omega$, $S(n)$ is a finite set
        of size  $n$. 
\item If $f$ is a function with $\dom(f) = \omega $, $S$ a
        slalom, then we say that $S$ {\em captures} $f$ iff
        $\forall^\infty n \,\, f(n) \in S(n)$. 
\item Let $M \subseteq  N$ be sets (typically: models of ZFC*).
We say that $N$ has the {\em Laver property} over $M$ iff:
\begin{quote} For every function $H\in \pre \omega^\omega
\cap M$, for every function $f\in \pre \omega^\omega \cap N$
satisfying $f \le H$ there is a slalom $S\in M$ that
captures $f$.
\end{quote}
\item A forcing notion $P$ has the Laver property iff
$\forces_P $ ``$V^P$ has the Laver property over $V$.''
\end{enumerate}
\end{Definition}



Before we formulate the main lemma, we  need the
following easy claim: 

\begin{Claim}
  Let $k> 2^n$. If $X \subseteq \pre  2^k$, 
        $|X| = n$ then there are $i<j$ in $k$ 
        such that for all $f\in X$, $f(i)=f(j)$. 
\end{Claim}

\begin{proof}
For $i<j$, $f\in X$,  define an equivalence relation $\sim_f$ by:
$i \sim_f j \iff   f(i) = f(j)$.
Let $i \sim j$ iff $i \sim_f j$ for all $f$ in $X$.   Since each $\sim_f$
has at most $2$ equivalence classes, $\sim$ has at most $2^n$ classes,
so there are $i \not= j$, $i \sim j$. 
\end{proof}

\begin{Lemma}\label{lavermain}
Assume that $(N,{\in} )$ is a model of ZFC*,  and that 
$V$ has the Laver property over $N$.

Then every real avoids some $2$-rule from  $N$. 
\end{Lemma}

\begin{proof}

Let $a_0 = 0$, $a_{n+1} = a_n + 2^n+1$.  The sequence 
$(a_n:n\in \omega )$ is in $N$. 

For any $X\in {\Cal P}(\omega)$, we will find a rule in $N$ which 
$X$ does not follow.

Let $\chi_X\in \ccc$ be the characteristic function of $X$. 
Define $X^* := (\chi_X \on [a_n, a_{n+1}): n \in \omega)$.  
Note that there are only $2^{2^n+1}$ many possibilities for
$\chi_X\on [a_n, a_{n+1} )$.

Since $V$ has the Laver
property over $N$ there is a sequence $\vec S = 
(S_n: n\in \omega)\in N$, 
$S_n \subseteq \pre 2^{[a_n,a_{n+1} )}$, 
$|S_n| \le n$, and for all 
$n>0$,  $\chi_X\on [a_n, a_{n+1}) \in S(n)$.     
By the above claim we can find $i_n< j_n$ in 
$[a_n, a_{n+1}) $ such that for all $z\in S(n)$, 
$z(i_n)=z(j_n)$.  Since  the sequence $\vec S$ is in $N$,
we can find such a sequence 
$( i_n, j_n:n<\omega)$ in $N$. 

Define a 2-rule $( A_n, B_n: n \in \omega)\in N$ by $A_n=\{i_n\}$, 
$B_n = \{i_n, j_n \}$.    Since $i_n\in X $ iff $j_n \in X$, 
$X $
does not follow this rule. 

\end{proof}

\begin{Lemma}\label{laverfact}
\begin{enumerate}
\itm a Let $\bar P = (P_i, Q_i: i< \omega_2 )$ be a countable support
        iteration of proper forcing notions such that for
        each $i$ we have $\forces_i$ ``$Q_i$ has the Laver
        property.''  Then $P_{\omega_2}$, the countable support
        limit of $\bar P$, also has the Laver property.
\itm b Laver forcing  is proper and  has the Laver property.
\itm c Laver forcing adds a real that dominates all reals
       from the ground model. 
\end{enumerate}
\end{Lemma}

\begin{proof}  These facts are 
  well known and (at least for the case where each $Q_i$
        is Laver forcing) appear implicitly or explicitly
in Laver's  paper \cite{L1}.
\end{proof}

\begin{Conclusion}
Let $P_{\omega_2} $ be the limit of a countable support
iteration of Laver forcing over a model $V_0$ of GCH.
Then $\forces_{P_{\omega_2}} \b=\r={\omega_2} $ and
$\r_2={\omega_1}$.
\end{Conclusion}
\begin{proof}
Let $V_{{\omega_2}} = V^{P_{\omega_2}}$. 
$V_{\omega_2} \models \b={\omega_2} $ is 
well known. (Let $f_i$ be the real added by
the $i$th Laver forcing; then $( f_i:i<{\omega_2})$ is a
strictly increasing and cofinal sequence in $\pre
\omega^\omega$.

By \ref{laverfact}, $V_{{\omega_2}}$ has the Laver property
over $V_0$.   Hence, by \ref{lavermain}, every real avoids
some rule from $V_0$.
So $\r_2 \le |\pre \omega^\omega \cap V_0|  = \aleph_1 $.
\end{proof}

%
%

%
%

\section{Application to independent families}

A family $\F\su\Cal P(\om)$ of subsets of $\om$ is {\em independent} iff
it generates a free boolean algebra in $\Cal P(\om)/\fin$. Equivalently,
for any two disjoint finite subsets of $\F$, the intersection of all
members in the first set with all complements of members in the second
set is infinite. 

The following
is an example of an independent family of size continuum over a
countable set: $\{A_r:r\in \Bbb R\}$ where $A_r=\{p\in \Bbb Z[X]:
p(r)>0\}$.

A family  $\F\su \Cal P(\om)$ is {\em dense} iff for any two
finite disjoint subsets of $\om$ there are infinitely many members of
$\F$ that contain the first set and are disjoint to the second.

An interesting (proper) subclass of the class of dense independent
families over $\om$ is the class of homogeneous families, which was
introduced in \cite{GGK}. Its study was continued in \cite{KjSh:499}.

While every dense independent family is contained in a maximal dense
independent family, this is not obvious (and perhaps false)  for homogeneous
families. The existence, even the consistency, of a maximal homogeneous
family over $\om$ is still open. In particular, an increasing union of
homogeneous families need not be homogeneous. 

In the study of extendibility of homogeneous families, the following
notion is fundamental: 
Let $G\su
\aut \F$. We define $( \F,G)\le ( \F',G')$ iff $\F\su \F'$,
$G\su G'\su \aut \F'$. The usefulness of $\le$ is that unions of
suitable $\le$ chains {\em are} homogeneous (see \cite{KjSh:499} for a
detailed account of direct limits in the category of homogeneous
families). 

We show now that below $\r_\infty$ one can get proper $\le$-extensions of
independent families. This was our original motivation for discovering
$\r_\infty$. 

\begin{theorem}\label{motivate}
Suppose $G \subseteq \aut \F$, $\F\su\Cal P(\om)$ is dense independent and
$|\F|+|G|<\r_\infty $.  Then there exists $\F'\supsetneqq \F$ such
that $( \F,G) \le ( \F',G)$
\end{theorem}

\begin{proof}
Suppose that $G\su \aut \F$, $\F$ is dense independent and
$|G|+|\F|<r_\infty$.
We shall find a real $X\su \om$ such that $X\notin \F$ and $\F\cup
G[X]$ is independent, where $G[X]$ is the orbit of $X$ under $G$. This
will suffice, since clearly $G\su \aut( \F\cup G[X])$ for any real $X$.

It is a priori unclear why such $X$ should exist. If for example there
is some $\sigma\in G$ with finite support, then for no $X\su \om$ is even
the orbit $G[X]$ itself independent. However, the following lemma
takes care of this.
Let $\supp({\sigma} ) = \{ n\in \omega : {\sigma}(n)\not=n\}$ for a
permutation ${\sigma} \in \sym \omega $.

\begin{lemma}\label{support}  Suppose that $\F$ is dense independent
and $\sigma\in \aut \F$ is not the identity. Then there are distinct
sets $C_n\in \F$ such that for all $n$, $C_{2n}-C_{2n+1}\su \supp
\sigma$. 
\end{lemma}

\begin{Remark}In particular, the supports of non-identity automorphisms
have the finite intersection property and hence generate a filter.
 This is the ``strong Mekler condition'' for
$\aut \F$ (see  \cite{truss}).
\end{Remark}

\begin{proof}
Fix $k^*\in \om$ for which $\sigma(k^*)\not=k^*$. Find $C_{2n}$, $C_{2n+1}$
by induction on $n$. 
Suppose $C_m$ is chosen for
$m< 2n$. By density, there are infinitely many $C\in \F$ for which
$k^*\in C$, $ {\sigma} (k^*)\notin C$. Choose some such $C$ so that 
neither $C$ nor $ {\sigma} [C]$
are  among $\{C_m:m< 2n\}$. Let $C_{2n}=C$ and
$C_{2n+1}= {\sigma} [C]$. Since $k^*\in C_{2n}-C_{2n+1}$, those sets are
indeed distinct.

We claim that $C_{2n}-C_{2n+1} \subseteq \supp( {\sigma} )$.  Indeed,
for any $k\in C_{2n}-C_{2n+1}$ we have $ {\sigma} (k)\in 
C_{2n+1}$ but $k \notin C_{2n+1}$, so $k \not= {\sigma} (k)$. 
\end{proof}
We now continue the proof of theorem \ref{motivate}. 
Let $M$ be a transitive model of ZFC* of cardinality $< \r_\infty$
such that 
$\F,G\in M$ and $G\su M, \F\su M$. Let $X$ be a real that follows
all bounded rules from $M$.  Clearly, $X\notin M$, and therefore $X\notin
\F$.

We need to show that every boolean combination over $\F\cup G[X]$
is infinite.
Suppose that
$$D:= A\cap \sigma_0[X] \cap\cdots\cap \sigma_{n-1}[X]\cap
(\om-\sigma_{n}[X])\cap\cdots \cap(\om-\sigma_{m-1}[X])$$
is a boolean combination over $\F\cup G[X]$, where $\sigma_i\in G$ for
$i<m$, and $A$ is some boolean combination over $\F$. Clearly, $A\in
M$.

Set $N=  \binom{m}{2}$ and let $( \tau_i:i<N)$ be a list of all
$\sigma_k\circ \sigma^{-1}_\ell$ for $k<\ell<m$.
By induction find a sequence $(C_0, \ldots, C_{2n+1})$ of
$2N+2$ many distinct sets such that no $C_k$ participates in $D$
and such that 
$C_{2k}-C_{2k+1}\su \supp \tau_k$.
$C_{2k}$ and $C_{2k+1}$ are constructed in the $k$-th step by using 
 lemma \ref{support}.

Since all the $C_k$ are distinct, 
the intersection $E=A\cap \bigcap_{k<N}C_{2k}-\bigcup_{k<N}C_{2k+1}$ is
infinite.   Clearly $E$ belongs to $M$.
Define by induction an $m$-rule $( A_n,B_n:n<\om)$ as
follows: suppose $( A_k,B_k:k\le n)$ are defined. Find a point
$j_n\in E$ such that $B= \{\sigma_k^{-1}(j_n): k<m\}$ is disjoint from
$\bigcup_{\ell\le n} B_\ell$ and $j_n \notin
\bigcup_{\ell\le n} B_\ell$.
Let $B_{n+1}$ be $B$ and let
$A_{n+1}=\{\sigma_\ell^{-1}(j_n): \ell< n\}$.

The rule we defined obviously belongs to $M$. Since $X$ satisfies all
bounded rules from $M$, there are infinitely many $n$ for which $X\cap
B_n=A_n$. For each such $n$, $X_n \in D$.
\end{proof}

\end{document}